\documentclass[a4paper,english,fontsize=10pt,parskip=half,abstracton]{scrartcl}
\usepackage{babel}
\usepackage[utf8]{inputenc}
\usepackage[T1]{fontenc}
\usepackage[a4paper,left=20mm,right=20mm,top=30mm,bottom=30mm]{geometry}
\usepackage{amsmath}
\usepackage{amsthm}
\usepackage{amssymb}
\usepackage{enumerate}
\usepackage[bookmarks=false,
            pdftitle={On a theorem of Blichfeldt},
            pdfauthor={Benjamin Sambale},
            pdfkeywords={Blichfeldt's Theorem, number of fixed points, permutation character},
            pdfstartview={FitH}]{hyperref}

\newtheorem{Theorem}{Theorem} 

\numberwithin{equation}{section}

\allowdisplaybreaks[1]

\renewcommand{\phi}{\varphi}

\newcommand{\Aff}{\operatorname{Aff}}

\newcommand{\SL}{\operatorname{SL}}

\mathchardef\ordinarycolon\mathcode`\:  
 \mathcode`\:=\string"8000
 \begingroup \catcode`\:=\active
   \gdef:{\mathrel{\mathop\ordinarycolon}}
 \endgroup

\title{On a theorem of Blichfeldt}
\author{Benjamin Sambale\footnote{Fachbereich Mathematik, TU Kaiserslautern, 67653 Kaiserslautern, Germany, 
\href{mailto:sambale@mathematik.uni-kl.de}{sambale@mathematik.uni-kl.de}}}
\date{\today}

\begin{document}
\frenchspacing
\maketitle
\begin{abstract}\noindent
Let $G$ be a permutation group on $n<\infty$ objects. Let $f(g)$ be the number of fixed points of $g\in G$, and let $\{f(g):1\ne g\in G\}=\{f_1,\ldots,f_r\}$. In this expository note we give a character-free proof of a theorem of Blichfeldt which asserts that the order of $G$ divides $(n-f_1)\ldots(n-f_r)$. We also discuss the sharpness of this bound.
\end{abstract}


Let us consider a permutation group $G$ on a finite set $\Omega$ consisting of $n$ elements. By Lagrange's Theorem applied to the symmetric group on $\Omega$, it follows that the order $|G|$ of $G$ is a divisor of $n!$. In order to strengthen this divisibility relation we denote the number of fixed points of a subgroup $H\le G$ on $\Omega$ by $f(H)$. Moreover, let $f(g):=f(\langle g\rangle)$ for every $g\in G$. In 1895, Maillet~\cite{Maillet} proved the following (see also Cameron's book \cite[p.~172]{Cameron}).

\begin{Theorem}[Maillet]
Let $\{f(H):1\ne H\le G\}=\{f_1,\ldots,f_r\}$. Then $|G|$ divides $(n-f_1)\ldots(n-f_r)$.
\end{Theorem}

Using the newly established character theory of finite groups, Blichfeldt~\cite{Blichfeldt} showed in 1904 that it suffices to consider cyclic subgroups $H$ in Maillet's Theorem (this was rediscovered by Kiyota~\cite{KiyotaBlich}).

\begin{Theorem}[Blichfeldt]\label{blich}
Let $\{f(g):1\ne g\in G\}=\{f_1,\ldots,f_r\}$. Then $|G|$ divides $(n-f_1)\ldots(n-f_r)$.
\end{Theorem}

For the convenience of the reader we present the elegant argument which can be found in \cite[Theorem~6.5]{Cameron}.

\begin{proof}[Proof of Blichfeldt's Theorem]
Since $f$ is the permutation character, the function $\psi$ sending $g\in G$ to $(f(g)-f_1)\ldots(f(g)-f_r)$ is a generalized character of $G$ (i.\,e. a difference of ordinary complex characters). From
\[\psi(g)=\begin{cases}
(n-f_1)\ldots(n-f_r)&\text{if }g=1,\\0&\text{if }g\ne 1
\end{cases}\]
we conclude that $\psi$ is a multiple of the regular character $\rho$ of $G$. In particular, $\rho(1)=|G|$ divides $\psi(1)=(n-f_1)\ldots(n-f_r)$.
\end{proof}

It seems that no elementary proof (avoiding character theory) of Blichfeldt's Theorem has been published so far. The aim of this note is to provide such a proof.

\begin{proof}[Character-free proof of Blichfeldt's Theorem]
It suffices to show that
\[\frac{1}{|G|}\sum_{g\in G}(f(g)-f_1)\ldots(f(g)-f_r)\in\mathbb{Z},\]
since all summands with $g\ne 1$ vanish. Expanding the product we see that it is enough to prove
\[F_k(G):=\frac{1}{|G|}\sum_{g\in G}{f(g)^k}\in\mathbb{Z}\]
for $k\ge 0$. Obviously, $F_0(G)=1$. Arguing by induction on $k$ we may assume that $F_{k-1}(H)\in\mathbb{Z}$ for all $H\le G$. Let $\Delta_1,\ldots,\Delta_s$ be the orbits of $G$ on $\Omega$, and let $\omega_i\in\Delta_i$ for $i=1,\ldots,s$. For $\omega\in\Delta_i$ the stabilizers $G_\omega$ and $G_{\omega_i}$ are conjugate in $G$. In particular, $F_{k-1}(G_\omega)=F_{k-1}(G_{\omega_i})$. Recall that the orbit stabilizer theorem gives us $|\Delta_i|=|G:G_{\omega_i}|$ for $i=1,\ldots,s$. This implies
\begin{align*}
F_k(G)&=\frac{1}{|G|}\sum_{\omega\in\Omega}\sum_{g\in G_\omega}{f(g)^{k-1}}=\frac{1}{|G|}\sum_{\omega\in\Omega}|G_\omega|F_{k-1}(G_\omega)=\frac{1}{|G|}\sum_{i=1}^s\sum_{\omega\in\Delta_i}|G_\omega|F_{k-1}(G_\omega)\\
&=\frac{1}{|G|}\sum_{i=1}^s|\Delta_i||G_{\omega_i}|F_{k-1}(G_{\omega_i})=\frac{1}{|G|}\sum_{i=1}^s|G:G_{\omega_i}||G_{\omega_i}|F_{k-1}(G_{\omega_i})=\sum_{i=1}^sF_{k-1}(G_{\omega_i})\in\mathbb{Z}.\qedhere
\end{align*}
\end{proof}

As a byproduct of the proof we observe that $F_1(G)$ is the number of orbits of $G$. This is a well-known formula sometimes (inaccurately) called Burnside's Lemma (see \cite{Burnsidenot}). If there is only one orbit, the group is called \emph{transitive}. In this case, $F_2(G)$ is the \emph{rank} of $G$, i.\,e. the number of orbits of any one-point stabilizer.

It is known that Blichfeldt's Theorem can be improved by considering only the fixed point numbers of non-trivial elements of prime power order. This can be seen as follows. Let $S_p$ be a Sylow $p$-subgroup of $G$ for every prime divisor $p$ of $|G|$. Since 
\[\{f(g):1\ne g\in S_p\}\subseteq\{f(g):1\ne g\in G\text{ has prime power order}\}=\{f_1,\ldots,f_r\},\] 
Theorem~\ref{blich} implies that $|S_p|$ divides $(n-f_1)\ldots(n-f_r)$ for every $p$. Since the orders $|S_p|$ are pairwise coprime, also $|G|=\prod_p|S_p|$ is a divisor of $(n-f_1)\ldots(n-f_r)$. On the other hand, it does not suffice to take the fixed point numbers of the elements of prime order. An example is given by $G=\langle(1,2)(3,4),(1,3)(2,4),(1,2)(5,6)\rangle$. This is a dihedral group of order $8$ where every involution moves exactly four letters.

Cameron-Kiyota~\cite{CamKiy} (and independently Chillag~\cite{Chillag}) obtained another generalization of Theorem~\ref{blich} where $f$ is assumed to be any generalized character $\chi$ of $G$ and $n$ is replaced by its degree $\chi(1)$. 
A dual version for conjugacy classes instead of characters appeared in Chillag~\cite{Chillag2}.

Numerous articles addressed the question of equality in Blichfeldt's Theorem. Easy examples are given by the \emph{regular} permutation groups. These are the transitive groups whose order coincides with the degree. In fact, by Cayley's Theorem every finite group is a regular permutation group acting on itself by multiplication. A wider class of examples consists of the \emph{sharply $k$-transitive} permutation groups $G$ for $1\le k\le n$. Here, for every pair of tuples $(\alpha_1,\ldots,\alpha_k),(\beta_1,\ldots,\beta_k)\in\Omega^k$ with $\alpha_i\ne\alpha_j$ and $\beta_i\ne\beta_j$ for all $i\ne j$ there exists a unique $g\in G$ such that $\alpha_i^g=\beta_i$ for $i=1,\ldots,k$. Setting $\alpha_i=\beta_i$ for all $i$, we see that any non-trivial element of $G$ fixes less than $k$ points. Hence,
\[\{f(g):1\ne g\in G\}\subseteq\{0,1,\ldots,k-1\}.\]
On the other hand, if $(\alpha_1,\ldots,\alpha_k)$ is fixed, then there are precisely $n(n-1)\ldots(n-k+1)$ choices for $(\beta_1,\ldots,\beta_k)$. It follows that $|G|=n(n-1)\ldots(n-k+1)$. Therefore, we have equality in Theorem~\ref{blich}. Note that sharply $1$-transitive and regular are the same thing. 
An interesting family of sharply $2$-transitive groups comes from the \emph{affine groups}
\[\Aff(1,p^m)=\{\phi:\mathbb{F}_{p^m}\to\mathbb{F}_{p^m}\mid\exists a\in\mathbb{F}_{p^m}^\times,b\in\mathbb{F}_{p^m}:\phi(x)=ax+b\ \forall x\in\mathbb{F}_{p^m}\}\]
where $\mathbb{F}_{p^m}$ is the field with $p^m$ elements.
More generally, all sharply $2$-transitive groups are \emph{Frobenius groups} with abelian kernel. By definition, a Frobenius group $G$ is transitive and satisfies $\{f(g):1\ne g\in G\}=\{0,1\}$. 
The \emph{kernel} $K$ of $G$ is the subset of fixed point free elements together with the identity. Frobenius Theorem asserts that $K$ is a (normal) subgroup of $G$. For the sharply $2$-transitive groups this can be proved in an elementary fashion (see \cite[Exercise~1.16]{Cameron}), but so far no character-free proof of the full claim is known.
The dihedral group $\langle(1,2,3,4,5),(2,5)(3,4)\rangle$ of order $10$ illustrates that not every Frobenius group is sharply $2$-transitive.

A typical example of a sharply $3$-transitive group is $\SL(2,2^m)$ with its natural action on the set of one-dimensional subspaces of $\mathbb{F}_{2^m}^2$. We leave this claim as an exercise for the interested the reader.
The sharply $k$-transitive groups for $k\in\{2,3\}$ were eventually classified by Zassenhaus~\cite{Zassenhaus1,Zassenhaus2} using near fields (see Passman's book~\cite[Theorems~20.3 and 20.5]{Passmanbook}). 
On the other hand, there are not many sharply $k$-transitive groups when $k$ is large. In fact, there is a classical theorem by Jordan~\cite{Jordan} which was supplemented by Mathieu~\cite{Mathieu}.

\begin{Theorem}[Jordan, Mathieu]
The sharply $k$-transitive permutation groups with $k\ge 4$ are given as follows:
\begin{enumerate}[(i)]
\item the symmetric group of degree $n\ge 4$ \textup{(}$k\in\{n,n-1\}$\textup{)},
\item the alternating group of degree $n\ge 6$ \textup{(}$k=n-2$\textup{)},
\item the Mathieu group of degree $11$ \textup{(}$k=4$\textup{)},
\item the Mathieu group of degree $12$ \textup{(}$k=5$\textup{)}.
\end{enumerate}
\end{Theorem}

We remark that the Mathieu groups of degree $11$ and $12$ are the smallest members of the \emph{sporadic} simple groups.

In accordance with these examples, permutation groups with equality in Theorem~\ref{blich} are now called \emph{sharp} permutation groups (this was coined by Ito-Kiyota~\cite{ItoK}).
Apart from the ones we have already seen, there are more examples. For instance, the symmetry group of a square acting on the four vertices has order $8$ (again a dihedral group) and the non-trivial fixed point numbers are $0$ and $2$. Recently, Brozovic~\cite{Brozovic} gave a description of the primitive sharp permutation groups $G$ such that $\{f(g):1\ne g\in G\}=\{0,k\}$ for some $k\ge 1$. Here, a permutation group is \emph{primitive} if it is transitive and any one-point stabilizer is a maximal subgroup. The complete classification of the sharp permutation groups is widely open.

Finally, we use the opportunity to mention a related result by Bochert~\cite{Bochert} where the divisibility relation of $|G|$ is replaced by an inequality. As usual $\lfloor x\rfloor$ denotes the largest integer less than or equal to $x\in\mathbb{R}$.

\begin{Theorem}[Bochert]
If $G$ is primitive, then $|G|\le n(n-1)\ldots(n-\lfloor n/2\rfloor+1)$ unless $G$ is the symmetric group or the alternating group of degree $n$.
\end{Theorem}

\section*{Acknowledgment}
This work is supported by the German Research Foundation (project SA 2864/1-1) and the Daimler and Benz Foundation (project 32-08/13).


\begin{thebibliography}{10}

\bibitem{Blichfeldt}
H.~F. Blichfeldt, \textit{A theorem concerning the invariants of linear
  homogeneous groups, with some applications to substitution-groups}, Trans.
  Amer. Math. Soc. \textbf{5} (1904), 461--466.

\bibitem{Bochert}
A. Bochert, \textit{Ueber die {Z}ahl der verschiedenen {W}erthe, die eine
  {F}unction gegebener {B}uchstaben durch {V}ertauschung derselben erlangen
  kann}, Math. Ann. \textbf{33} (1889), 584--590.

\bibitem{Brozovic}
D.~P. Brozovic, \textit{The classification of primitive sharp permutation
  groups of type {$\{0,k\}$}}, Comm. Algebra \textbf{42} (2014), 3028--3062.

\bibitem{Cameron}
P.~J. Cameron, \textit{Permutation groups}, London Mathematical Society Student
  Texts, Vol. 45, Cambridge University Press, Cambridge, 1999.

\bibitem{CamKiy}
P.~J. Cameron and M. Kiyota, \textit{Sharp characters of finite groups}, J.
  Algebra \textbf{115} (1988), 125--143.

\bibitem{Chillag}
D. Chillag, \textit{Character values of finite groups as eigenvalues of
  nonnegative integer matrices}, Proc. Amer. Math. Soc. \textbf{97} (1986),
  565--567.

\bibitem{Chillag2}
D. Chillag, \textit{On a congruence of {B}lichfeldt concerning the order of
  finite groups}, Proc. Amer. Math. Soc. \textbf{136} (2008), 1961--1966.

\bibitem{ItoK}
T. Ito and M. Kiyota, \textit{Sharp permutation groups}, J. Math. Soc. Japan
  \textbf{33} (1981), 435--444.

\bibitem{Jordan}
C. Jordan, \textit{Sur la limite de transitivité des groupes non-alternées},
  Bull. Soc. Math. France \textbf{1} (1873), 40--71.

\bibitem{KiyotaBlich}
M. Kiyota, \textit{An inequality for finite permutation groups}, J. Combin.
  Theory Ser. A \textbf{27} (1979), 119.

\bibitem{Maillet}
E. Maillet, \textit{Sur quelques propri\'et\'es des groupes de substitutions
  d'ordre donn\'e}, Ann. Fac. Sci. Toulouse Sci. Math. Sci. Phys. \textbf{9}
  (1895), 1--22.

\bibitem{Mathieu}
{\'E}. Mathieu, \textit{Mémoire sur l'étude des fonctions de plusieurs
  quantités, sur la mani\`{e}re de les former et sur les substitutions qui les
  laissent invariables}, J. Math. Pures Appl. \textbf{6} (1861), 241--323.

\bibitem{Burnsidenot}
P.~M. Neumann, \textit{A lemma that is not {B}urnside's}, Math. Sci. \textbf{4}
  (1979), 133--141.

\bibitem{Passmanbook}
D.~S. Passman, \textit{Permutation groups}, Dover Publications, Inc., Mineola,
  N.Y., 2012 (revised reprint of the 1968 original).

\bibitem{Zassenhaus2}
H. Zassenhaus, \textit{Kennzeichnung endlicher linearer {G}ruppen als
  {P}ermutationsgruppen}, Abh. Math. Sem. Univ. Hamburg \textbf{11} (1935),
  17--40.

\bibitem{Zassenhaus1}
H. Zassenhaus, \textit{\"{U}ber endliche {F}astk\"orper}, Abh. Math. Sem. Univ.
  Hamburg \textbf{11} (1935), 187--220.

\end{thebibliography}
\end{document}